\documentclass[11pt,a4paper,reqno]{amsart}
\usepackage[T1]{fontenc}
\usepackage[utf8]{inputenc}
\usepackage{lmodern}
\usepackage{amsmath,amssymb,mathtools}
\mathtoolsset{showonlyrefs}
\usepackage{microtype}
\usepackage[margin=30mm]{geometry}
\usepackage[hidelinks]{hyperref}
\hypersetup{
  pdftitle={Radial solutions of pure Neumann Henon-Lane-Emden systems on balls and the weighted critical hyperbola},
  pdfauthor={Qiaoran Wu and Wenming Zou},
  pdfsubject={Weighted radial nonexistence and a qualitative blow-up consequence},
  pdfkeywords={Henon-Lane-Emden system, pure Neumann boundary conditions, weighted critical hyperbola, radial nonexistence}
}
\newtheorem{theorem}{Theorem}[section]
\newtheorem{lemma}[theorem]{Lemma}

\newtheorem{corollary}[theorem]{Corollary}
\theoremstyle{remark}
\newtheorem{remark}[theorem]{Remark}
\numberwithin{equation}{section}
\newcommand{\R}{\mathbb R}

\newcommand{\clBR}{\overline{B_R}}
\newcommand{\dd}{\,\mathrm d}
\newcommand{\norm}[1]{\left\lVert #1\right\rVert}
\newcommand{\doi}[1]{\href{https://doi.org/#1}{\nolinkurl{#1}}}
\allowdisplaybreaks[1]
\title[Radial solutions and the weighted critical hyperbola]{{Radial solutions of pure Neumann
H\'enon--Lane--Emden systems on balls and the weighted critical hyperbola}}

\author{Qiaoran Wu}
\address{Department of Mathematical Sciences, Tsinghua University, Beijing, China}
\email{wuqr24@mails.tsinghua.edu.cn}

\author{Wenming Zou}
\address{Department of Mathematical Sciences, Tsinghua University, Beijing, China}
\email{zou-wm@mail.tsinghua.edu.cn}
\date{}

\begin{document}

\begin{abstract}
{
We study radial solutions of the pure Neumann H\'enon--Lane--Emden system
$-\Delta u=|x|^\alpha|v|^{q-1}v$, $-\Delta v=|x|^\beta|u|^{p-1}u$
in a ball $B_R\subset\R^N$, where $N\geq3$, $p,q>0$, and $\alpha,\beta>-2$.
Solutions are assumed to be continuous on the closed ball and classical
away from the centre. We prove that every radial solution is trivial if
\[
\frac{N+\alpha}{q+1}+\frac{N+\beta}{p+1}\leq N-2.
\]
In the unweighted case, this resolves the question of Pistoia, Schiera,
and Tavares concerning radial solutions on the critical hyperbola in a ball.
When $pq\ne1$, the strict reverse inequality is necessary and sufficient
for the existence of a nontrivial radial solution. In this subcritical
regime, the radial least-energy pair is unique up to simultaneous sign
reversal and is radially nondegenerate. After a choice of sign, both radial
profiles have strictly negative derivatives on $(0,R)$, and each has
exactly one simple zero. The nonexistence proof combines a weighted
radial Rellich--Pohozaev identity with a local uniqueness argument for
the radial system at criticality. For fixed weights, if a sequence of
subcritical exponent pairs converges to a positive point of the weighted
critical hyperbola, the $L^\infty$ norms of both components of any
corresponding sequence of nontrivial radial solutions diverge.
We also derive a nonexistence result for a weighted fourth-order Neumann problem.
\par}
\end{abstract}

\maketitle

\section{Introduction and main result}

Let $N\geq3$ and $R>0$, and let $B_R=B_R(0)\subset\R^N$.
We study the H\'enon--Lane--Emden system
\begin{equation}\label{eq:system}
\begin{cases}
-\Delta u=|x|^\alpha|v|^{q-1}v &\text{in }B_R\setminus\{0\},\\
-\Delta v=|x|^\beta|u|^{p-1}u &\text{in }B_R\setminus\{0\},\\
\partial_\nu u=\partial_\nu v=0 &\text{on }\partial B_R,
\end{cases}
\qquad p,q>0,
\end{equation}
where $\nu$ denotes the outward unit normal. Unless explicitly stated
otherwise, the weights satisfy $\alpha,\beta>-2$. A solution of \eqref{eq:system} means
a pair $u,v\in C(\clBR)\cap C^2(\clBR\setminus\{0\})$ satisfying the
equations away from the origin and the Neumann conditions on $\partial B_R$.
For nonnegative weights, every radial solution in this class is a
classical $C^2(\clBR)$ solution; see Remark~\ref{rem:originregularity}.
A solution is
nontrivial if $(u,v)\ne(0,0)$, and is radial if both components are
radial with respect to the centre of the ball.

The weighted critical hyperbola for the radial problem is
\begin{equation}\label{eq:critical}
\frac{N+\alpha}{q+1}+\frac{N+\beta}{p+1}=N-2.
\end{equation}
This curve also appears in the study of radial H\'enon--Lane--Emden
systems on the punctured whole space; see Musina and Sreenadh
\cite{MS14}, whose notation for the exponents differs from ours.
We call the regime
$(N+\alpha)/(q+1)+(N+\beta)/(p+1)>N-2$ subcritical and the regime
$(N+\alpha)/(q+1)+(N+\beta)/(p+1)<N-2$ supercritical.

\begin{samepage}
Integrating on annuli and passing to the limit
at the origin gives
\begin{equation}\label{eq:compatibility}
\int_{B_R}|x|^\beta|u|^{p-1}u\dd x
=\int_{B_R}|x|^\alpha|v|^{q-1}v\dd x=0.
\end{equation}
\end{samepage}
If either component vanishes identically, so does the other. Hence both
components of every nontrivial radial solution change sign.

In the unweighted case $\alpha=\beta=0$, the subcritical condition is
$1/(p+1)+1/(q+1)>(N-2)/N$. In this regime, Salda\~na and Tavares established existence,
regularity, and qualitative properties of least-energy nodal solutions
when $pq\ne1$ \cite[Theorem~1.1]{ST18}. On balls and annuli, they also
proved symmetry breaking for least-energy solutions in the full class:
these solutions are foliated Schwarz symmetric but not radial. At the
same time, minimization restricted to the radial class yields classical
radial solutions, which are strictly monotone in the radial variable
\cite[Theorem~1.2]{ST18}.

On the unweighted critical hyperbola, Pistoia, Schiera, and Tavares proved existence
of least-energy nodal classical solutions under dimension-dependent
restrictions on the exponents, using a dual variational formulation and
critical compactness estimates \cite[Theorems~1.2 and~1.5]{PST24}. Under
the corresponding hypotheses on balls and annuli, they also proved that
least-energy solutions are nonradial \cite[Theorem~1.9(ii)]{PST24}.
This symmetry-breaking result concerns only the least-energy level.
For $p\ne q$, it leaves open whether a ball admits any radial
solution at criticality \cite[Remark~1.10]{PST24};
see also \cite[Section~3.2, p.~190]{T24}, which records it for least-energy
radial solutions. The
diagonal critical point, $p=q=(N+2)/(N-2)$, is already understood: the
system reduces to the scalar critical Neumann equation, whose Pohozaev
argument rules out radial solutions \cite[Remark~6.5]{PST24}. That
argument does not extend off the diagonal, since it does not ensure that
both components vanish on the same concentric sphere. Our theorem resolves
the off-diagonal case when $\alpha=\beta=0$ and proves radial
nonexistence throughout the weighted critical and supercritical regimes.
In the unweighted slightly supercritical regime on the unit
ball, sign-changing solutions do exist \cite{GP25}; by
Theorem~\ref{thm:main}, none of them is radial.

Our main result is the following.

\begin{theorem}\label{thm:main}
Let $N\geq3$, $R>0$, $p,q>0$, and $\alpha,\beta>-2$. Assume that
\begin{equation}\label{eq:nonexistcondition}
\frac{N+\alpha}{q+1}+\frac{N+\beta}{p+1}\leq N-2.
\end{equation}
Then every radial solution $(u,v)$ of \eqref{eq:system} with
$u,v\in C(\clBR)\cap C^2(\clBR\setminus\{0\})$ is trivial.
Consequently, neither component of a nontrivial solution in this class
can be radial.
\end{theorem}

{
The following theorem characterizes the exponents for which a nontrivial
radial solution exists and describes the least-energy pair in the radial class.

\begin{theorem}\label{thm:subcritical}
Let $N\geq3$, $R>0$, $\alpha,\beta>-2$, and $p,q>0$ with $pq\ne1$.
Then \eqref{eq:system} admits a nontrivial radial solution in the stated
class if and only if
\[
  \frac{N+\alpha}{q+1}+\frac{N+\beta}{p+1}>N-2.
\]
Under this condition, there exists a radial least-energy solution.
This least-energy pair is unique up to simultaneous sign reversal. After
a choice of sign, both profiles have strictly negative derivatives
on $(0,R)$. Each has one simple nodal radius, hence exactly two nodal
domains. The integrated radial linearization at this pair has trivial
Neumann kernel (see \eqref{eq:sub-linear}).
\end{theorem}

For $\alpha=\beta=0$, the criterion reduces to
$1/(p+1)+1/(q+1)>(N-2)/N$. In this case the existence direction also
follows from Salda\~na and Tavares \cite[Theorem~1.2]{ST18}.
The uniqueness of the radial least-energy pair answers the ball case of
the question raised at the end of their introduction.
\par}

The restriction $pq\ne1$ in \cite[Theorem~1.2]{ST18} is genuine:
$pq=1$ is the asymptotically linear (hence homogeneous) case.
The formal reduction leads to a higher-order equation
(see \cite[eq.~(3.4), p.~182]{T24}). For Neumann boundary conditions, the corresponding auxiliary
nonlinear eigenvalue problem is studied in \cite[Section~2]{SST26}, which
also analyzes the limiting behavior of least-energy solutions as $pq\to1$
\cite[Theorems~1.1--1.2]{SST26}. This curve is strictly subcritical
and is not covered by Theorem~\ref{thm:main}.

For $\alpha=\beta=0$ at criticality, the main theorem also rules out radial solutions in the
strong solution class considered by Pistoia, Schiera and Tavares
\cite{PST24}.

\begin{corollary}\label{cor:strong}
Let $N\geq4$, $\alpha=\beta=0$, and let $p,q>0$ satisfy
$1/(p+1)+1/(q+1)=(N-2)/N$.
Suppose that
\[
(u,v)\in W^{2,(q+1)/q}(B_R)\times W^{2,(p+1)/p}(B_R)
\]
solves \eqref{eq:system}. If either component is radial, then
$u=v=0$ in $B_R$.
\end{corollary}

Section~\ref{sec:consequences} also records a fourth-order consequence
(Corollary~\ref{cor:biharmonic}): for $N\geq5$, $\beta>-2$, and
$p\geq(N+4+2\beta)/(N-4)$ the Neumann problem
$\Delta^2w=|x|^\beta|w|^{p-1}w$ in $B_R\setminus\{0\}$, with
$\partial_\nu w=\partial_\nu(\Delta w)=0$ on $\partial B_R$, has no
nontrivial radial solution in $C^4(\clBR)$.

For the unweighted scalar pure Neumann equation the near-critical picture is
finer: sharp radial concentration estimates are available
\cite{GST20}, and least-energy solutions converge to a least-energy
solution of the critical problem \cite{ST22}. For radial solutions
of the weighted system approaching the weighted critical hyperbola from
the subcritical side, we obtain the following qualitative blow-up result.

\begin{theorem}\label{thm:blowup}
Let $N\geq3$, $R>0$, and fix $\alpha,\beta>-2$.
Let $(p_n,q_n)_{n\geq1}$ be a sequence of positive exponent pairs satisfying
\[
\frac{N+\alpha}{q_n+1}+\frac{N+\beta}{p_n+1}>N-2,
\]
such that
\[
(p_n,q_n)\longrightarrow(p_*,q_*)\in(0,\infty)^2,
\qquad
\frac{N+\alpha}{q_*+1}+\frac{N+\beta}{p_*+1}=N-2.
\]
For each $n$, let $(u_n,v_n)$ be a nontrivial radial solution of
\eqref{eq:system}, with
$u_n,v_n\in C(\clBR)\cap C^2(\clBR\setminus\{0\})$, exponents $(p_n,q_n)$,
and the fixed weights $(\alpha,\beta)$. Then
\[
\norm{u_n}_\infty\longrightarrow\infty,
\qquad
\norm{v_n}_\infty\longrightarrow\infty.
\]
\end{theorem}

{ Theorem~\ref{thm:subcritical} supplies nontrivial radial
solutions along every such subcritical sequence for which $p_nq_n\ne1$.
The blow-up statement applies to any nontrivial radial solutions, not
just to the least-energy pairs.\par}

The argument uses the classical Rellich--Pohozaev method for elliptic
systems; see Mitidieri \cite{M93} and van der Vorst \cite{V91}.
In the present setting the relevant identity is
\begin{equation}\label{eq:introidentity}
\begin{aligned}
&\frac{R^{N+\beta}|U(R)|^{p+1}}{p+1}
 +\frac{R^{N+\alpha}|V(R)|^{q+1}}{q+1}\\
&\qquad=\left(\frac{N+\beta}{p+1}
              +\frac{N+\alpha}{q+1}-(N-2)\right)A,
\end{aligned}
\end{equation}
where $u(x)=U(|x|)$, $v(x)=V(|x|)$ and
\[
\begin{aligned}
A&=\int_0^R r^{N-1}U'(r)V'(r)\dd r\\
 &=\int_0^R r^{N-1+\beta}|U(r)|^{p+1}\dd r
  =\int_0^R r^{N-1+\alpha}|V(r)|^{q+1}\dd r.
\end{aligned}
\]
Thus $A>0$ for a nontrivial solution. The supercritical case gives an
immediate sign contradiction. At criticality,
\eqref{eq:introidentity} forces both boundary values to vanish, and the
Neumann conditions supply $U'(R)=V'(R)=0$ as well. Uniqueness for the
radial ODE system with these zero Cauchy data then implies that the
solution is trivial.

For nonnegative weight exponents, the weighted Dirichlet problem in a
ball carries the same exponent condition: Moreira dos Santos, Nornberg and
Soave prove that, under
\eqref{eq:nonexistcondition}, no nontrivial radial classical solution
vanishing on $\partial B_R$ exists \cite[Theorem~1.8]{MNS21}. Up to
interchanging the two weights, which those authors attach to the opposite
equations, their identity carries the same coefficient and the same integral
as \eqref{eq:introidentity}, and differs only on the left. For Dirichlet
data that side is $R^NU'(R)V'(R)$, a quantity the identity does not sign, so
positivity has to be recovered from a unique continuation principle
\cite[Corollary~2.3]{MNS21}. Under the Neumann condition it is the boundary
energy
$R^{N+\beta}|U(R)|^{p+1}/(p+1)+R^{N+\alpha}|V(R)|^{q+1}/(q+1)$,
whose nonnegativity is immediate. Thus no separate argument is needed to
determine the sign of the boundary term. At criticality, the proof still
uses the elementary uniqueness result for zero Cauchy data in
Lemma~\ref{lem:zero}.

Section~\ref{sec:proof} proves Theorem~\ref{thm:main}.
Section~\ref{sec:consequences} proves Corollary~\ref{cor:strong} in the unweighted case, gives the weighted
fourth-order application, and discusses the role of the ball. Section~\ref{sec:blowup} proves
Theorem~\ref{thm:blowup} and compares its unweighted specialization with
the convergence results for unrestricted least-energy solutions in
\cite{SST26}.
{ Section~\ref{sec:subcritical} proves existence, uniqueness
and radial nondegeneracy in the weighted subcritical regime.\par}

\section{Proof of radial nonexistence}\label{sec:proof}

For $t>0$ write $\phi_t(s)=|s|^{t-1}s$, with $\phi_t(0)=0$.
The radial equations associated with \eqref{eq:system} are
\begin{equation}\label{eq:radial}
\begin{cases}
-U''-\dfrac{N-1}{r}U'=r^\alpha\phi_q(V),\\[1mm]
-V''-\dfrac{N-1}{r}V'=r^\beta\phi_p(U),
\end{cases}
\qquad 0<r<R.
\end{equation}
The boundary conditions are $U'(R)=V'(R)=0$. No derivative condition
is imposed at the origin.

\begin{lemma}\label{lem:zero}
Let $p,q>0$ with $pq>1$ and let $\alpha,\beta\in\R$. Suppose that
$U,V\in C^2((0,R])\cap C([0,R])$ solve \eqref{eq:radial} and that
\[
U(R)=V(R)=U'(R)=V'(R)=0.
\]
Then $U\equiv V\equiv0$ on $[0,R]$.
\end{lemma}

\begin{proof}
We first establish local uniqueness of the zero solution at any
$r_0\in(0,R]$. Suppose that
$U(r_0)=V(r_0)=U'(r_0)=V'(r_0)=0$.
For $0<\delta<r_0/2$, put
\[
J_\delta=[r_0-\delta,\min\{r_0+\delta,R\}],\qquad
M_U(\delta)=\max_{J_\delta}|U|,\quad
M_V(\delta)=\max_{J_\delta}|V|.
\]
Integrating \eqref{eq:radial} from $r_0$ yields
\begin{equation}\label{eq:localintegral}
U'(r)=-r^{1-N}\int_{r_0}^{r}t^{N-1+\alpha}\phi_q(V(t))\dd t,
\qquad r\in J_\delta,
\end{equation}
and the analogous formula for $V'$, with $\beta$ and $p$ in place of
$\alpha$ and $q$.
Since $r,t\in[r_0/2,3r_0/2]$, the powers of $r$ and $t$ are bounded
there. A second integration gives a constant
$C>0$, independent of small $\delta$, such that
\begin{equation}\label{eq:localestimates}
M_U(\delta)\leq C\delta^2M_V(\delta)^q,
\qquad
M_V(\delta)\leq C\delta^2M_U(\delta)^p.
\end{equation}
If $M_U(\delta)>0$, these inequalities imply
\begin{equation}\label{eq:zerocontradiction}
1\leq C^{1+q}\delta^{2(1+q)}M_U(\delta)^{pq-1}.
\end{equation}
The right-hand side tends to zero as $\delta\downarrow0$, because
$pq>1$ and $U$ is continuous at $r_0$. Consequently,
$M_U(\delta)=0$ for all sufficiently small $\delta$; then
\eqref{eq:localestimates} also gives $M_V(\delta)=0$.

The set
\[
Z=\{r\in(0,R]:U(r)=V(r)=U'(r)=V'(r)=0\}
\]
is nonempty, is relatively closed by continuity, and is relatively open
by the preceding local argument. Since $(0,R]$ is connected,
$Z=(0,R]$. Continuity at the origin completes the proof.
(A radial unique continuation theorem for much wider classes, including
$pq\ge1$ and Pucci operators, is \cite[Theorem~1.7]{MNS21}. The elementary
argument above is used below with singular weights, where that theorem does
not apply.)
\end{proof}

\begin{lemma}\label{lem:identity}
Let $N\geq3$, $R>0$, $p,q>0$, and $\alpha,\beta>-2$.
Let $(u,v)$ be radial, continuous on $\clBR$ and of class $C^2$ in
$\clBR\setminus\{0\}$, and satisfy \eqref{eq:system} pointwise away from
the origin. Define
\[
A=\int_0^Rr^{N-1}U'V'\dd r,
\]
where $u(x)=U(|x|)$ and $v(x)=V(|x|)$. Then
\begin{equation}\label{eq:cross}
A=\int_0^Rr^{N-1+\beta}|U|^{p+1}\dd r
 =\int_0^Rr^{N-1+\alpha}|V|^{q+1}\dd r,
\end{equation}
and \eqref{eq:introidentity} holds.
\end{lemma}

\begin{proof}
Since $U,V$ are bounded and $N+\alpha,N+\beta>0$, the derivatives of
$r^{N-1}U'$ and $r^{N-1}V'$ are integrable near zero. Therefore, $r^{N-1}U'$ and $r^{N-1}V'$ have finite limits at zero. A nonzero limit would make $U$ and $V$ unbounded because $N\geq3$, so both limits vanish. Integrating the equations \eqref{eq:radial} gives
    \begin{equation}\label{eq:originintegral}
    	U'(r)=-r^{1-N}\int_0^r t^{N-1+\alpha}\phi_q(V(t))\dd t
    	=O(r^{\alpha+1}),
    \end{equation}
and
    \begin{equation}\label{eq:V-deriviation}
    	V'(r)=-r^{1-N}\int_0^r t^{N-1+\beta}\phi_p(U(t))\dd t
    	=O(r^{\beta+1}).
    \end{equation}
At $r=R$, the Neumann conditions in these formulas also give
\eqref{eq:compatibility}.
Thus $r^{N-1}U'V$ and $r^{N-1}V'U$ tend to zero as $r\downarrow0$.
Moreover, $r^{N-1}|U'V'|=O(r^{N+\alpha+\beta+1})$ is integrable near zero
because $N+\alpha+\beta+1>-1$. The two weighted integrals in
\eqref{eq:cross} are also finite, since $U,V$ are bounded and
$N+\alpha,N+\beta>0$.

Multiplying the first equation in \eqref{eq:radial} by
$r^{N-1}V$ and integrating by parts on $(\varepsilon,R)$ gives
\[
\int_\varepsilon^Rr^{N-1+\alpha}|V|^{q+1}\dd r
=-[r^{N-1}U'V]_\varepsilon^R
 +\int_\varepsilon^Rr^{N-1}U'V'\dd r.
\]
Letting $\varepsilon\downarrow0$ and using the Neumann condition at
$R$, we obtain the equality with the last integral in \eqref{eq:cross}.
Testing the second equation with $r^{N-1}U$ on $(\varepsilon,R)$ and passing
to the same limit gives the remaining equality.

Set
\[
F(s)=\frac{|s|^{p+1}}{p+1},\qquad
G(s)=\frac{|s|^{q+1}}{q+1},
\]
\[
H(r)=U'(r)V'(r)+r^\beta F(U(r))+r^\alpha G(V(r)).
\]
The functions $F,G$ are continuously differentiable for every $p,q>0$.
Using \eqref{eq:radial}, we obtain, for $r>0$,
\[
H'(r)=-\frac{2(N-1)}rU'(r)V'(r)
       +\beta r^{\beta-1}F(U(r))
       +\alpha r^{\alpha-1}G(V(r)).
\]
Hence
\begin{equation}\label{eq:differentialidentity}
\begin{aligned}
\frac{\mathrm d}{\mathrm dr}\bigl(r^NH(r)\bigr)
={}&(2-N)r^{N-1}U'V'\\
 &+(N+\beta)r^{N-1+\beta}F(U)
  +(N+\alpha)r^{N-1+\alpha}G(V).
\end{aligned}
\end{equation}
The same estimates give
$r^NH(r)=O(r^{N+\alpha+\beta+2})+O(r^{N+\beta})+O(r^{N+\alpha})\to0$, all
exponents being positive because $N\geq3$ and $\alpha,\beta>-2$.
Integrating \eqref{eq:differentialidentity} on $(\varepsilon,R)$ and
letting $\varepsilon\downarrow0$, then using $U'(R)=V'(R)=0$ and
\eqref{eq:cross}, we obtain \eqref{eq:introidentity}.
\end{proof}

\begin{remark}\label{rem:originregularity}
For $\alpha,\beta\geq0$, every radial solution in the class used above
is classical also at the origin. Indeed, set
$h(r)=r^\alpha\phi_q(V(r))$, with its continuous value at $r=0$. \eqref{eq:originintegral} gives, for $r>0$,
\[
\begin{aligned}
\frac{U'(r)}r&=-\int_0^1 t^{N-1}h(rt)\dd t,\\
U''(r)&=(N-1)\int_0^1 t^{N-1}h(rt)\dd t-h(r).
\end{aligned}
\]
Both expressions tend to $-h(0)/N$ as $r\downarrow0$. Since
$U'(r)=O(r)$, the radial extension $u(x)=U(|x|)$ is differentiable at the
origin with zero gradient, and its Hessian extends continuously there.
Thus $u\in C^2(\clBR)$ and $-\Delta u(0)=h(0)$. The same argument applies
to $v$.
\end{remark}

\begin{proof}[Proof of Theorem~\ref{thm:main}]
Suppose that $(u,v)$ is a nontrivial radial solution in the class of
Theorem~\ref{thm:main}. By \eqref{eq:cross}, $A>0$.
If the inequality in \eqref{eq:nonexistcondition} is strict, the
right-hand side of \eqref{eq:introidentity} is negative, whereas its left-hand side is
nonnegative, a contradiction.

If equality holds in \eqref{eq:nonexistcondition}, then
\eqref{eq:introidentity} gives $U(R)=V(R)=0$. Moreover,
\[
(N-2)\left(\frac1{p+1}+\frac1{q+1}\right)
<\frac{N+\alpha}{q+1}+\frac{N+\beta}{p+1}=N-2
\quad\Longrightarrow\quad pq>1.
\]
This verifies the hypothesis $pq>1$ required in
Lemma~\ref{lem:zero}.
Lemma~\ref{lem:zero}, together with the Neumann conditions, yields
$U=V=0$, again a contradiction.

Finally, if only $u$ is assumed radial, then
$\phi_q(v)=|x|^{-\alpha}(-\Delta u)$ is radial for $x\ne0$.
Since $\phi_q$ is one-to-one, $v$ is radial away from the origin and
hence on the whole ball by continuity. The same argument applies with
the components interchanged.
\end{proof}

\section{Consequences and the scope of the result}\label{sec:consequences}

\subsection{Unweighted critical strong solutions}

\begin{proof}[Proof of Corollary~\ref{cor:strong}]
Here $\alpha=\beta=0$. We adopt the setting of \cite{PST24}, where
$N\geq4$; the two spaces in
the corollary are those of their strong-solution class. By Proposition~1.3 of that paper, every such strong solution belongs
to $C^{2,\zeta}(\clBR)\times C^{2,\eta}(\clBR)$ for any
$0<\zeta<\min\{q,1\}$ and $0<\eta<\min\{p,1\}$, and this regularity
conclusion holds throughout the unweighted critical hyperbola. Theorem~\ref{thm:main}
therefore applies.
\end{proof}

\subsection{A weighted fourth-order consequence}

\begin{corollary}\label{cor:biharmonic}
Let $N\geq5$, $R>0$, $\beta>-2$, and
$p\geq (N+4+2\beta)/(N-4)$. Then the problem
\begin{equation}\label{eq:biharmonic}
\begin{cases}
\Delta^2w=|x|^\beta|w|^{p-1}w &\text{in }B_R\setminus\{0\},\\
\partial_\nu w=\partial_\nu(\Delta w)=0 &\text{on }\partial B_R
\end{cases}
\end{equation}
has no nontrivial radial solution in $C^4(\clBR)$.
\end{corollary}

\begin{proof}
Given such a radial solution $w$, define $u=w$ and $v=-\Delta w$.
Then $(u,v)$ is a radial solution of \eqref{eq:system} in the class of
Theorem~\ref{thm:main}, with $\alpha=0$ and $q=1$. The assumed lower bound on $p$ is equivalent to
\[
\frac{N+\beta}{p+1}+\frac N2\leq N-2.
\]
Theorem~\ref{thm:main} now gives the claimed nonexistence.
\end{proof}

\section{Proof of radial blow-up near criticality}\label{sec:blowup}

To prove Theorem~\ref{thm:blowup}, we first derive an amplitude estimate
that rules out convergence to the zero solution in the compactness
argument.
Throughout this section $\alpha,\beta>-2$ are fixed and all supremum
norms are taken over $B_R$.

\begin{lemma}\label{lem:amplitude}
Let $(u,v)$ be a nontrivial radial solution of \eqref{eq:system} in
the solution class specified above, with $p,q>0$. Put
\[
M_U=\norm{u}_\infty,\qquad M_V=\norm{v}_\infty,
\]
\[
C_\alpha=\frac{R^{\alpha+2}}{(N+\alpha)(\alpha+2)},\qquad
C_\beta=\frac{R^{\beta+2}}{(N+\beta)(\beta+2)}.
\]
Then
\begin{equation}\label{eq:amplitude}
M_U\leq C_\alpha M_V^q,\qquad M_V\leq C_\beta M_U^p.
\end{equation}
In particular, if $pq>1$, then
\begin{equation}\label{eq:lowerbounds}
M_U\geq(C_\alpha C_\beta^q)^{-1/(pq-1)},\qquad
M_V\geq(C_\beta C_\alpha^p)^{-1/(pq-1)}.
\end{equation}
\end{lemma}

\begin{proof}
By \eqref{eq:originintegral} and \eqref{eq:V-deriviation}, we have
\[
|U'(r)|\leq\frac{r^{\alpha+1}}{N+\alpha}M_V^q,
\qquad |V'(r)|\leq\frac{r^{\beta+1}}{N+\beta}M_U^p.
\]
By \eqref{eq:compatibility}, the radial functions $U$ and $V$ each have
at least one zero in $(0,R)$.
Since $\alpha+2,\beta+2>0$, the derivative bounds are integrable at
zero. Integrating the bound for $U'$ from a zero of $U$ and the bound for
$V'$ from a zero of $V$ proves \eqref{eq:amplitude}; these zeros need not
coincide.
Combining the inequalities gives
\[
M_U\leq C_\alpha C_\beta^q M_U^{pq},\qquad
M_V\leq C_\beta C_\alpha^p M_V^{pq}.
\]
Since $M_U,M_V>0$, \eqref{eq:lowerbounds} follows.
\end{proof}

\begin{proof}[Proof of Theorem~\ref{thm:blowup}]
Because $N+\alpha,N+\beta>N-2$, the weighted critical condition
implies $1/(p_*+1)+1/(q_*+1)<1$, hence $p_*q_*>1$.
Since $p_nq_n\to p_*q_*>1$, there is $\delta>0$ such that
$p_nq_n\geq1+\delta$ for all sufficiently large $n$.
Suppose, to the contrary, that $\norm{u_n}_\infty$ does not tend
to infinity. Then there are $M>0$ and a subsequence, still indexed by
$n$, such that $\norm{u_n}_\infty\leq M$. Since $p_n,q_n\to p_*,q_*>0$,
the exponents $p_n,q_n$ remain in a compact subinterval of $(0,\infty)$, so
the numbers $M^{p_n}$ are uniformly bounded. The second inequality in
\eqref{eq:amplitude} therefore gives a uniform upper bound for
$\norm{v_n}_\infty$ along this subsequence. Since $C_\alpha,C_\beta$
are fixed positive constants and the exponents converge with
$p_nq_n\geq1+\delta$, the lower bounds \eqref{eq:lowerbounds} are
uniform in $n$ and give positive bounds for both component norms.

Let $U_n,V_n$ be the radial functions associated with $u_n,v_n$.
Applying \eqref{eq:originintegral}
and its analogue for $V_n'$, the uniform amplitude bounds give, for
$0\leq s\leq r\leq R$,
\[
\begin{aligned}
|U_n(r)-U_n(s)|&\leq C\bigl(r^{\alpha+2}-s^{\alpha+2}\bigr),\\
|V_n(r)-V_n(s)|&\leq C\bigl(r^{\beta+2}-s^{\beta+2}\bigr),
\end{aligned}
\]
where $C$ is independent of $n$. Since both powers are positive,
these estimates give equicontinuity on $[0,R]$. By the
Arzel\`a--Ascoli theorem, after taking a further subsequence,
\[
U_n\longrightarrow U,\qquad V_n\longrightarrow V
\quad\text{uniformly on }[0,R].
\]
The uniform lower bounds on the component norms ensure that neither
$U$ nor $V$ is identically zero. By the compactness of the exponents
established above, together with the uniform bounds on $U_n,V_n$, the pairs
$(q_n,V_n(r))$ and $(p_n,U_n(r))$ lie in a common compact subset of
$(0,\infty)\times\R$, on which $(t,s)\mapsto\phi_t(s)$ is uniformly
continuous. Consequently,
\[
\phi_{q_n}(V_n)\longrightarrow\phi_{q_*}(V),\qquad
\phi_{p_n}(U_n)\longrightarrow\phi_{p_*}(U)
\quad\text{uniformly on }[0,R].
\]
For $r>0$, rewrite \eqref{eq:originintegral} as
\begin{equation}\label{eq:rescaledintegral}
U_n'(r)=-r^{\alpha+1}\int_0^1t^{N-1+\alpha}
                         \phi_{q_n}(V_n(rt))\dd t,
\end{equation}
and similarly for $V_n'$, with $\beta$ and $p_n$ in place of
$\alpha$ and $q_n$. Since $N+\alpha,N+\beta>0$, the powers of $t$ in
these integrals are integrable. The uniform convergence of the nonlinear
terms therefore gives uniform convergence of the derivatives on
$[\varepsilon,R]$ for every $\varepsilon\in(0,R)$. Hence
$U,V\in C^1((0,R])$ and satisfy the limiting integral formulas for
$U'$ and $V'$. Differentiating these formulas for $r>0$ shows that
$U,V\in C^2((0,R])$
and solve \eqref{eq:radial} with exponents $(p_*,q_*)$.
The derivative convergence at $R$ also gives $U'(R)=V'(R)=0$.

The radial extensions of $U,V$ are therefore continuous on $\clBR$,
of class $C^2$ away from the origin, and form a nontrivial solution of the
critical system. This contradicts Theorem~\ref{thm:main} and proves
$\norm{u_n}_\infty\to\infty$. Interchanging the components proves
$\norm{v_n}_\infty\to\infty$.
\end{proof}

\begin{remark}
In the unweighted case, sharper radial concentration estimates are known for the scalar pure
Neumann problem \cite{GST20}, while its least-energy solutions converge to
a least-energy solution of the critical problem \cite{ST22}. Theorem~\ref{thm:blowup}
applies to all sequences of radial solutions satisfying its hypotheses,
also for fixed power weights $\alpha,\beta>-2$.
\end{remark}

\begin{remark}
In the unweighted case, the least-energy solutions considered in
\cite{SST26} are nonradial under
the symmetry hypotheses of \cite[Theorem~1.1]{ST18} and
\cite[Theorem~1.9(ii)]{PST24}, so their convergence is consistent with
Theorem~\ref{thm:blowup}. Our theorem excludes uniformly $L^\infty$-bounded
radial sequences approaching a positive critical point from the
subcritical side. The radial least-energy sequences with $p_n\to0$ and
$q=1$ in \cite[Section~5, p.~990]{SST26} have a degenerate limit outside
our assumptions.
\end{remark}

\section{Proof of the radial existence criterion and least-energy properties}\label{sec:subcritical}

With $\phi_t$ as in Section~\ref{sec:proof}, set
\[
 \begin{gathered}
 \omega(r)=r^{N-1},\qquad d\mu_\gamma(r)=r^{N-1+\gamma}\dd r,\\
 \sigma=\frac{N+\beta}{p+1},\qquad \tau=\frac{N+\alpha}{q+1},
 \qquad \delta=\sigma+\tau-(N-2).
 \end{gathered}
\]
Here $\delta>0$. Throughout this section we assume $pq\ne1$.

\begin{lemma}\label{lem:sub-existence}
There exists a radial least-energy pair. Up to simultaneous sign reversal,
both profiles of every such pair have strictly negative derivatives on
$(0,R)$, and each profile has exactly one simple nodal radius.
\end{lemma}

\begin{proof}

Let $s=(p+1)/p$, $t=(q+1)/q$, and let $X_\beta$ and
$X_\alpha$ be the zero-mean subspaces of $L^s(\mu_\beta)$
and $L^t(\mu_\alpha)$, respectively. For $f\in X_\beta$ and
$g\in X_\alpha$, write
\[
 F(r)=\int_0^r f\,d\mu_\beta,\qquad
 G(r)=\int_0^r g\,d\mu_\alpha,\qquad
 \mathcal{P}(f,g)=\int_0^R \frac{F(r)G(r)}{\omega(r)}\dd r .
\]
Maximize $\mathcal{P}$ on the product of the closed unit balls of
$X_\beta$ and $X_\alpha$. H\"older's inequality gives
$|F(r)|\le C r^{\sigma}$ and $|G(r)|\le C r^{\tau}$ on these balls. Thus
the integrand is bounded by $C r^{\delta-1}$. Weak convergence
gives pointwise convergence of the primitives (the integration
interval can be tested by its characteristic function), so dominated
convergence and weak compactness give a maximizer. The maximum
$m$ is positive: one can take both primitives equal to a
nonzero smooth function compactly supported in $(0,R)$ and then
normalize the sources. By bilinearity both maximizing norms are one.

For a maximizer put $u=\phi_{1/p}(f)$ and $v=\phi_{1/q}(g)$.
Differentiating the separate unit-sphere constraints in their
zero-mean spaces gives
\[
 \mathcal{P}(h,g)=m\int_0^R uh\,d\mu_\beta\quad(h\in X_\beta),
 \qquad
 \mathcal{P}(f,k)=m\int_0^R vk\,d\mu_\alpha\quad(k\in X_\alpha).
\]
Indeed, the derivative of the norm at a norm-one source $f$
is integration against $\phi_{s-1}(f)=u$; the multiplier is
found by taking $h=f$, and likewise for $g$.
For $\zeta\in C_c^\infty(0,R)$ take
$h=\zeta'/r^{N-1+\beta}$, whose primitive is $\zeta$. Using the
analogous test for $g$ yields, first in distributions and then
pointwise for continuous representatives on $(0,R]$,
\begin{equation}\label{eq:sub-euler}
 u'=-\frac{G}{m\omega},\qquad v'=-\frac{F}{m\omega}.
\end{equation}
The primitives vanish at both endpoints. Continuity in the interior
implies that $f=\phi_p(u)$ and $g=\phi_q(v)$ have continuous
representatives there, and the profiles are $C^2$ up to $R$ away
from zero. Set $U=c u$, $V=d v$, where
$c=m d^q$, $d=m c^p$. These equations have a unique positive
solution because $pq\ne1$. The rescaled profiles solve the radial
equations with zero outer derivatives and zero weighted derivatives $\omega U'$, $\omega V'$ at the origin.

We verify continuity at the origin, including for negative weights.
From \eqref{eq:sub-euler} and the bounds on the primitives, integration
from a fixed positive radius gives $U=O(r^{-\sigma})$ and
$V=O(r^{-\tau})$, since $N-2-\tau=\sigma-\delta$ and $N-2-\sigma=\tau-\delta$.
Fix $0<\varepsilon<\min\{\delta,\alpha+2,\beta+2\}$. If the same
bounds hold with $(a,b)$ in place of $(\sigma,\tau)$, for
$0\le a\le \sigma$, $0\le b\le \tau$, the integrated equations give
\[
 |U'(r)|\le C r^{1+\alpha-qb},\qquad
 |V'(r)|\le C r^{1+\beta-pa}.
\]
Indeed $N+\alpha-qb\ge \tau>0$ and $N+\beta-pa\ge \sigma>0$.
Integrating these bounds gives new growth exponents
\[
 a'=\max\{0,qb-(\alpha+2)+\varepsilon\},\qquad
 b'=\max\{0,pa-(\beta+2)+\varepsilon\}.
\]
The extra $\varepsilon$ absorbs a possible logarithm. Starting at
$(\sigma,\tau)$, the iterates decrease coordinatewise: the first iterate
is at most $(\sigma,\tau)$ by $\varepsilon<\delta$, and the map is
order preserving. The limit is a fixed point. A fixed point with one
zero component is $(0,0)$. At a positive fixed point,
\[
 (pq-1)a=\alpha+2+q(\beta+2)-(q+1)\varepsilon
         =(pq-1)\sigma+(q+1)(\delta-\varepsilon).
\]
The first right-hand expression is positive, excluding $pq<1$;
the second then gives $a>\sigma$ when $pq>1$. Hence the iterates
tend to zero. Both become zero after finitely many steps, because
the two constants subtracted before taking the positive part are
strictly positive. Finally $U'=O(r^{1+\alpha})$ and
$V'=O(r^{1+\beta})$ are integrable at zero. The continuous
extensions therefore belong to the stated solution class.

To identify the least-energy pair, use the radial Hamiltonian energy
(with the area of the unit sphere suppressed)
\begin{equation}\label{eq:sub-energy}
 E(U,V)=\int_0^R U'V'\omega\dd r
 -\frac1{p+1}\int_0^R |U|^{p+1}d\mu_\beta
 -\frac1{q+1}\int_0^R |V|^{q+1}d\mu_\alpha .
\end{equation}
For any nontrivial radial solution, integration by parts gives a
common positive value $I$ for the three integrals (without their
coefficients). The boundary terms vanish by the zero weighted derivatives and the
bounded profiles. The solution's two source functions $\phi_p(U)$ and
$\phi_q(V)$ have zero means. Their primitives are $-\omega V'$
and $-\omega U'$, respectively. Dividing the sources by
$I^{p/(p+1)}$ and $I^{q/(q+1)}$ produces a norm-one pair
whose pairing is $I^{-\theta}$, where
$\theta=1-1/(p+1)-1/(q+1)\ne0$. The energy is $\theta I$.
As a function of $I^{-\theta}$ this energy is strictly decreasing:
its derivative with respect to $I^{-\theta}$ is
$-I^{1+\theta}$. The rescaling of a maximizing pair leaves its
normalized sources unchanged, so it attains the radial least energy.
Conversely, every radial least-energy pair gives a maximizer after
normalization.

The replacements $(F,G)\mapsto(|F|,|G|)$ preserve the constraints:
the new sources are the derivatives of the absolute primitives
divided by their weighted densities, and the magnitudes of the
derivatives are unchanged almost everywhere. The pairing cannot
decrease, so the new pair also maximizes. By \eqref{eq:sub-euler}
its two recovered profiles are nonincreasing. Their weighted source
means vanish, and each source is nonzero, so both profiles change
sign. The source of each nonincreasing profile is positive near zero
and negative near $R$. The primitive of that source is consequently
strictly positive throughout $(0,R)$: before the sign change this
follows by integrating from zero, and afterwards by integrating the
negative tail back from the zero value at $R$. Thus $|F|,|G|$
are strictly positive there. The original continuous primitives cannot
change sign in $(0,R)$, and their signs must agree because their
pairing is positive. Up to simultaneous sign reversal, both original
profiles have strictly negative derivatives throughout $(0,R)$.
In particular, each has exactly one simple nodal radius.
\end{proof}

\begin{lemma}\label{lem:sub-unique}
On a fixed ball there is at most one radial solution pair whose two
profiles are strictly decreasing on $(0,R)$.
\end{lemma}

\begin{proof}

Take two pairs with both profiles decreasing. Each has strictly
negative interior derivatives: the zero means of the sources and monotonicity
give the same strictly positive primitives just described. The
profiles are strictly positive at zero and strictly negative at $R$.
For
\[
  a=\frac{\alpha+2+q(\beta+2)}{pq-1},\qquad
  b=\frac{\beta+2+p(\alpha+2)}{pq-1},
\]
the transformation $(U,V)(r)\mapsto
(\lambda^a U(\lambda r),\lambda^b V(\lambda r))$ maps solutions
on radius $R$ to solutions on radius $R/\lambda$:
$a+2+\alpha=qb$ and $b+2+\beta=pa$. Since $a\ne0$, use
this scaling to set $U_1(0)=U_2(0)=1$. Denote the
new radii by $R_1,R_2$.

Suppose $V_2(0)>V_1(0)$. The first integrated equation \eqref{eq:originintegral} gives
$U_2'<U_1'$ while $V_2>V_1$, and hence $U_2<U_1$.
The second integrated equation \eqref{eq:V-deriviation} then gives $V_2'>V_1'$.
Thus the difference $V_2-V_1$ increases as long as it stays
positive, precluding a first contact. Both strict derivative
inequalities persist at $\min(R_1,R_2)$ by the integral formulas
and the strict ordering on the preceding interval. If $R_2\le R_1$,
then $0=U_2'(R_2)<U_1'(R_2)\le0$; if $R_1\le R_2$, then
$0=V_1'(R_1)<V_2'(R_1)\le0$. Both are contradictions. Exchanging
indices also excludes $V_1(0)>V_2(0)$.

For equal central data, the profiles coincide near the origin. The
two nonzero positive initial values keep the power functions locally
Lipschitz there, and their integral equations bound the supremum
norm of the difference on $[0,h]$ by
$C(h^{2+\alpha}+h^{2+\beta})$ times that norm. To continue equality
away from zero, for $0<s<1$ use
\[
 |\phi_s(x)-\phi_s(y)|\le C_s |x|^{s-1}|x-y|\qquad(x\ne0).
\]
For $|x-y|<|x|/2$ this follows by the mean value theorem. Otherwise,
the $s$-H\"older bound for $\phi_s$ gives it directly.
Taking $x$ from one reference profile gives locally integrable
difference coefficients at its simple zeros. When an exponent is at least one, the corresponding coefficient is bounded. The integrated first-order
system for the differences of profiles and weighted derivatives therefore admits
the Gr\"onwall inequality on compact intervals away from the origin.
The two normalized solutions coincide up to the smaller endpoint.
The radii must agree, or one profile would have zero derivative at
an interior point of the other. Undoing the scalings on the original
fixed ball gives the same pair.
\end{proof}

\begin{lemma}\label{lem:sub-kernel}
At a radial least-energy pair, the integrated radial linearization
defined in \eqref{eq:sub-linear} below has trivial Neumann kernel.
\end{lemma}

\begin{proof}

For the decreasing pair, define the radial linearization in the
integrated sense: its solution space consists of continuous radial
profiles $(X,Y)$ with
\begin{equation}\label{eq:sub-linear}
 \begin{aligned}
 \omega(r)X'(r)&=-\int_0^r q|V|^{q-1}Y\,d\mu_\alpha,\\
 \omega(r)Y'(r)&=-\int_0^r p|U|^{p-1}X\,d\mu_\beta .
 \end{aligned}
\end{equation}
The values of the coefficients at isolated zeros do not affect the
integrals. Kernel elements satisfy $X'(R)=Y'(R)=0$. The integral
formulation imposes zero weighted derivatives at the origin and accommodates exponents
below one since the simple interior zeros give integrable coefficients.
The linear solution space before imposing the outer conditions is
two-dimensional. Indeed, arbitrary prescribed initial values at zero
give a unique local solution by a small-interval contraction; away from zero the linear system for the profiles and weighted derivatives
has locally integrable coefficients and hence unique continuation for
every such initial vector.

For the solution with $(X(0),Y(0))=(0,1)$, the first equation of
\eqref{eq:sub-linear} gives $X'<0$ while $Y>0$, so $X<0$;
the second then gives $Y'>0$, precluding a first zero of $Y$.
At $R$ the same integral formulas yield $X'(R)<0<Y'(R)$.
A second solution is the scaling derivative
\[
 (X_d,Y_d)=(aU+rU',\ bV+rV').
\]
To see that it solves \eqref{eq:sub-linear}, differentiate the
scaled equations on intervals away from the zeros. The weighted derivatives of the
scaling derivative are continuous across the zeros, and their
derivatives in the linearized equations are integrable there; thus
the equations hold in integrated form across the zeros. At the
origin the weighted derivatives vanish: the radial equations and the estimates
$U'=O(r^{1+\alpha})$, $V'=O(r^{1+\beta})$ give
$X_d'=O(r^{1+\alpha})$, $Y_d'=O(r^{1+\beta})$.
The first initial value is $aU(0)\ne0$. The outer derivative
vector is $(R U''(R),R V''(R))$, because the unperturbed outer
derivatives vanish. Both entries are positive by the equations and
the negative boundary values of the profiles. This vector is linearly
independent of $(X'(R),Y'(R))$, whose entries have opposite signs.
Thus no nonzero linear combination satisfies both outer Neumann
conditions. The simultaneous sign reversal leaves the linearized
coefficients unchanged.
\end{proof}

\begin{proof}[Proof of Theorem~\ref{thm:subcritical}]
Theorem~\ref{thm:main} gives necessity of the strict subcritical
inequality. Conversely, Lemma~\ref{lem:sub-existence} gives existence
and the asserted nodal properties. After orienting each least-energy
pair so that both profiles decrease, Lemma~\ref{lem:sub-unique} gives
uniqueness. Radial nondegeneracy follows from Lemma~\ref{lem:sub-kernel}.
\end{proof}

\section*{Declarations}

\noindent\textbf{Funding.}
W.~Zou is supported by the National Key R\&D Program of China
(Grant No.~2023YFA1010001) and the National Natural Science Foundation of China
(Grant No.~12171265).

\medskip
\noindent\textbf{Conflict of interest.}
The authors declare that they have no conflicts of interest relevant to this
article.

\medskip
\noindent\textbf{Data availability.}
No datasets were generated or analysed in this theoretical study.

\medskip
\noindent\textbf{Use of generative artificial intelligence.}
The authors used ChatGPT (OpenAI) for literature searches, for the
development and checking of mathematical arguments, and for manuscript
drafting and editing. The authors take full responsibility for the
mathematical content, references, and final version.

\end{document}